\magnification=1200

\input amstex
\documentstyle{amsppt}

\vcorrection{-.4in}\hcorrection{.2 in}

\input xypic
\input epsf

\define\CDdashright#1#2{&\,\mathop{\dashrightarrow}\limits^{#1}_{#2}\,&}
\define\CDdashleft#1#2{&\,\mathop{\dashleftarrow}\limits^{#1}_{#2}\,&}

\def\P{\Bbb P}

\def\C{\Bbb C}

\def\Til#1{{\widetilde{#1}}}

\def\cmp{{c_{\text{SM}}}}
\def\cma{{c_{\text{Ma}}}}
\def\cwm{{c_{\text{wMa}}}}

\def\siltable#1.{
\vbox{\tabskip=0pt \offinterlineskip
\halign to360pt{\strut##& ##\tabskip=1em plus2em&
  \hfil##\hfil& \vrule##&
  \hfil##\hfil& \vrule##&
  \hfil##\hfil& \vrule##&
  \hfil##\hfil& \vrule\thinspace\vrule##&
  \hfil##\hfil& \vrule\thinspace\vrule##&
  \hfil##\hfil& ##\tabskip=0pt\cr
#1}}}


\CompileMatrices

\topmatter
\title Weighted Chern@-Mather classes and Milnor classes of hypersurfaces
\endtitle
\author Paolo Aluffi\endauthor
\date October 1998\enddate
\address 
Mathematics Department,
Florida State University,
Tallahassee, FL 32306
\endaddress
\email aluffi\@math.fsu.edu\endemail

\abstract 
We introduce a class extending the notion of Chern@-Mather class to
possibly nonreduced schemes, and use it to express the difference
between Schwartz@- MacPherson's Chern class and the class of the
virtual tangent bundle of a singular hypersurface of a nonsingular
variety.  Applications include constraints on the possible
singularities of a hypersurface and on contacts of nonsingular
hypersurfaces, and multiplicity computations.\endabstract
\endtopmatter

\document
\footnote[]{The author is grateful to Florida State University for a
`Developing Scholar Award' under which this research was completed}
\rightheadtext{Chern@-Mather classes and Milnor classes}

 
\head \S0. Introduction\endhead

The notion of Chern@-Mather class was first introduced in the
literature by Robert MacPherson, as one of the main ingredients in his
definition of functorial Chern classes for possibly singular complex
varieties. One way to think about Mather's class of $Y$ as defined by
MacPherson is the following: blow@-up $Y$ so that the pull@-back of
its sheaf of differentials is locally free modulo torsion; then mod
out the torsion, dualize, and take Chern classes. The operation can in
fact be performed for any sheaf; this is worked out in
\cite{Kwieci\'nski}.

This definition ignores possible nilpotents on $Y$. We feel that it
would be desirable to have a class in the spirit of Chern@-Mather
class, but in some way sensitive to possible nonreduced structures on
$Y$: first, this is natural from the algebro@-geometric standpoint;
secondly, as we will see, a natural candidate carries useful
information when applied to the {\it singularity subscheme\/} of a
hypersurface (for which possibly non@-reduced scheme structures play a
fundamental r\^ole).

Our candidate is introduced in \S1. Its definition is a suitable
weighted sum of `conventional' Chern@-Mather classes of subvarieties
of $Y$. The subvarieties are the supports of the components of the
(intrinsic) normal cone of $Y$, and the weights are the lengths of the
components of this cone. The class we obtain (trivially) agrees with
Mather's if $Y$ is a local complete intersection.

If $Y$ is the singularity subscheme of a hypersurface, we can relate
the weighted Chern@-Mather class with other natural classes defined in
this case. For example, in \cite{Aluffi1} we have defined and studied a
{\it $\mu$@-class\/} associated with the singularity subscheme of a
hypersurface; in this paper, we answer a question which we could not
address previously: how to give a reasonable definition for arbitrary
schemes $Y$, from which the $\mu$@-class could be recovered if $Y$ is
the singularity subscheme of a hypersurface $X$. The weighted
Chern@-Mather class is precisely such a class (Corollary~1.4). We hope
that this viewpoint will eventually give us the right hint on how to
define a $\mu$@-class for the singularities of more general varieties $X$.

The main application of weighted Chern@-Mather classes is to the
computation of the difference between MacPherson's class of the
hypersurface and the class of its virtual tangent bundle. A formula
for the difference, in terms of the $\mu$@-class, is proved
`numerically' in \cite{Aluffi3}, and at the level of Chow groups in
\cite{Aluffi2} (Theorem~I.5). Such differences have been named `Milnor
classes', as they generalize the fact that, for local complete
intersections with isolated singularities, the Milnor {\it number\/}
computes the difference between the (topological) Euler characteristic
and the degree of the class of the virtual tangent bundle (see
\cite{Suwa}, \cite{B@-L@-S@-S}). 

Weighted Chern@-Mather classes allow us to recast the formula from
\cite{Aluffi2}. We state this in \S1 (Theorem~1.2), together with
other facts regarding weighted Chern@-Mather classes, such as their
relation vis@-a@-vis a class appearing in \cite{P@-P} or their
behavior under blow@-ups. Proofs of these statements are sketched in
\S2, together with a few general considerations regarding Milnor
classes. Theorem~1.2 is proved in full in \S2.

The expression for the $\mu$@-class in terms of weighted Chern@-Mather
classes allows us in principle to compute the former for a wide class
of examples.  We give a couple of applications in this direction in
\S3, in the spirit of the examples worked out in \cite{Aluffi1}, \S4.
For example, we prove that if two nonsingular hypersurfaces $M_1$,
$M_2$ of degrees $d_1$, $d_2$ in projective space are tangent along a
positive dimensional subvariety, then $d_1=d_2$. This fact was proved
in \cite{Aluffi2}, but with a strong additional hypothesis on the
contact locus of $M_1$ and $M_2$; the new formula for the $\mu$@-class
shows that the extra hypothesis is unnecessary. We also collect in \S3
a few explicit computations of weighted Chern@-Mather classes.

\subheading{Acknowledgements} I am very grateful to Jean@-Paul
Brasselet and to Tatsuo Suwa for organizing the Sapporo symposium on 
`Singularities in Geometry and Topology'. Conversations with the
participants at the meeting, especially Piotr Pragacz and Shoji
Yokura, were very helpful.  I am particularly indebted to Piotr
Pragacz (and to the referee of \cite{Aluffi2}) for pointing out that
the main formula in \cite{Aluffi2} should be interpreted as the
computation of the characteristic cycle of a hypersurface. Adam
Parusi\'nski and Piotr Pragacz give an alternative proof of this
formula in \cite{P@-P}, by a local computation of multiplicities which
relates it to a formula from \cite{B@-M@-M} (over $\C$, and in
homology) for the characteristic cycle. The proof of Theorem~1.2
given in \S2 owes much to their approach: it is my attempt to produce
a proof in their style, but which does not rely on the complex
geometry of the situation (hence valid in the algebraic category and
for rational equivalence; and potentially more amenable to
generalizations, say to positive characteristic), or on the
sophisticated machinery of index formulas.

 
\head \S1. Weighted Chern@-Mather classes.\endhead

All schemes in this note are of finite type over an algebraically closed
field of characteristic~0, and (for simplicity) embeddable in an
ambient nonsingular variety, which we will denote by $M$.

Assume that $Y$ is reduced and irreducible, of dimension $k$.
The {\it Chern@-Mather\/} class of $Y$ can be defined as follows. Let
$G_k(TM)$ denote the Grassmann bundle whose fiber over $p\in M$
consists of the Grassmannian of $k$@-planes in $TM$, and let $Y^\circ$
be the nonsingular locus in $Y$. Consider the map
$$Y^\circ @>>> G_k(TM)$$
defined by sending $p\in Y^\circ$ to $T_pY\subset T_pM$. The {\it Nash
blow@-up\/} of $Y$ is the closure $\Til Y$ of the image of this map;
it comes equipped with a proper map $\nu$ to $Y$, and with the
restriction $T$ of the tautological subbundle over $G_k(TM)$. This
data is easily checked to be independent of the ambient variety $M$. 
The Chern@-Mather class of $Y$ is defined by
$$\cma(Y):=\nu_*\left(c(T)\cap [\Til Y]\right)$$
in the Chow group $A_*Y$ of $Y$. This class of course agrees with the
total (`homology') class of the tangent bundle of $Y$ if $Y$ happens
to be nonsingular to begin with.

Note that this definition assumes that $Y$ is reduced, as it needs $Y$
to be nonsingular at the general point, and ignores by construction
the presence of nilpotents along subvarieties of $Y$. Our task is to
modify this notion to take account of possible nilpotents on $Y$.

Let then $Y\subset M$ be arbitrary. We consider the normal cone $C_YM$
of $Y$ in $M$, and associate with $Y$ the set $\{(Y_i,m_i)\}_i$, where
the $Y_i\overset{j_i}\to\hookrightarrow Y$ are the supports of the
irreducible components $C_i$ of $C_YM$, and $m_i$ denotes the geometric
multiplicity of $C_i$ in $C_YM$ (so $[C_i]=m_i[(C_i)_{\text{red}}]$).
\proclaim{Lemma 1.1}
The data $\{(Y_i,m_i)\}$ is intrinsic of $Y$, i.e., independent of the
ambient nonsingular variety.\endproclaim

\demo{Proof} (Cf.~\cite{Fulton}, Example~4.2.6.) It is enough to compare
embeddings $Y\hookrightarrow M$, $Y\hookrightarrow M'$, where both
$M$, $M'$ are nonsingular, and $M$ is smooth over $M'$. In this case
there is an exact sequence of cones
$$0 @>>> T_{M'|M} @>>> C_YM @>>> C_YM' @>>> 0$$
(where $T_{M'|M}$ is the relative tangent bundle) in the sense of
\cite{Fulton}, Example~4.1.6, and it follows that the supports of the
irreducible components of the two cones coincide, as well as the
geometric multiplicities of the components.\qed\enddemo

By Lemma 1.1, the following definition is also intrinsic of $Y$:
\definition{Definition} The {\it weighted Chern@-Mather class\/} of
$Y$ is
$$\cwm(Y):=\sum_i (-1)^{\dim Y-\dim Y_i} m_i j_{i*}
\cma(Y_i)\quad\text{in $A_*Y$.}$$
(Warning: we will henceforth neglect to indicate `obvious' push@-forwards
such as $j_{i*}$, and pull@-backs.)\enddefinition

Note that if $Y$ is a reduced irreducible local complete intersection,
then its normal cone is reduced and irreducible, so the class defined
here agrees with the Chern@-Mather class of $Y$. In particular, if $Y$
is nonsingular then $\cwm(Y)=c(TY)\cap [Y]$ is the total homology
class of the tangent bundle of $Y$.\vskip 6pt

A few examples of computations of weighted Chern@-Mather classes can
be found in \S3. Our main motivation in introducing the class
$\cwm(Y)$ is that we can prove it is particularly well@-behaved if $Y$
is the {\it singularity scheme of a hypersurface\/} $X$ in a
nonsingular variety $M$. By hypersurface here we mean the zero@-scheme
of a nonzero section of a line@-bundle $\Cal L$ on $M$; the
singularity subscheme of $X$ is the subscheme locally defined by
the partial derivatives of an equation for $X$. (This scheme structure
is independent of the ambient variety $M$.) In the rest of this
section we survey a few facts about $\cwm(Y)$ under the hypothesis that
$Y$ is the singularity subscheme of a hypersurface. Proofs are given
in \S2.

Our motivation is to highlight apparently different contexts in which
the class $\cwm(Y)$ manifests itself. Although these contexts will
invoke other characters of the play, remember that $\cwm(Y)$ is a
class {\it intrinsic of $Y$,\/} and which is defined regardless of
whether $Y$ is the singularity subscheme of a hypersurface. The
challenge is to find extensions of these results which do not assume
that $Y$ is the singularity subscheme of a hypersurface.

For the first fact, let $\cmp(X)$ denote Schwartz@-MacPherson's Chern
class of $X$, and let $c_F(X)$ denote the class of its virtual tangent
bundle; the subscript $F$ is to remind us that this class agrees with
the class introduced (for much more general schemes) by William Fulton,
cf.~Example 4.2.6 of \cite{Fulton}.

\proclaim{Theorem 1.2} Let $\Cal L=\Cal O(X)$, and let $Y$ be the
singularity subscheme of $X$. Then
$$\cwm(Y)=(-1)^{\dim X-\dim Y} c(\Cal L)\cap \left(c_F(X)-\cmp(X)\right)
\quad\text{in $A_*(X)$.}$$
\endproclaim
That is, $\cwm(Y)$ essentially measures the difference between the
functorial homology Chern class $\cmp(X)$ and the class of the virtual
tangent bundle of $X$. The functoriality of the class $\cmp(X)$ was
proved by Robert MacPherson \cite{MacPherson}; the class was later
shown to agree with the class previously defined by Marie@-H\'el\`ene
Schwartz. For a treatment of Schwartz@-MacPherson's classes over any
algebraically closed field of characteristic~0, see \cite{Kennedy};
this is the context we assume here. Also, we let
$\cmp(X)=\cmp(X_{\text{red}})$; with this proviso, Theorem~1.2
holds for nonreduced hypersurfaces $X$---remarkably, the drastic
change in $c_F$ when some component of $X$ is replaced by a multiple
is precisely compensated by the change in the weighted Mather class of
the singularity subscheme.

For the next result, it is convenient to employ the following notations
(a variation on the notations used in \cite{Aluffi2}, \cite{Aluffi3}):
for $a\in A_p$ and $\Cal L$ a line bundle, set
$$a_\vee=(-1)^p a\quad,\quad a_{\Cal L}=c(\Cal L)^p\cap a\quad.$$
(So $a_{\Cal L}=c(\Cal L)^n\cap (a\otimes \Cal L)$, where the term in
$()$ uses the definition in \cite{Aluffi3}, and $n$ is the dimension of
the ambient scheme). These notations satisfy simple identities,
analogous to the ones shown in \cite{Aluffi3}. For example, the
formula on the right defines an action of Pic on the Chow group: that
is, $a_{\Cal L_1\otimes\Cal L_2}=(a_{\Cal L_1})_{\Cal L_2}$ for line
bundles $\Cal L_1$ and $\Cal L_2$.

\proclaim{Proposition 1.3}
$$\cwm(Y)=(-1)^{\dim Y}\left(c(T^*M\otimes \Cal L)\cap s(Y,M)\right)_{
\vee\Cal L}\quad\text{in $A_*Y$.}$$
\endproclaim
Here $s(Y,M)$ denotes the {\it Segre class\/} of $Y$ in $M$, in the
sense of \cite{Fulton}, Chapter~4. Note that this equality is
completely false unless $Y$ is a singularity subscheme of a
hypersurface in $M$. However, if $Y$ {\it is\/} a singularity subscheme
of a hypersurface in $M$, then the right@-hand@-side must be independent
of $M$: this was proved directly in \cite{Aluffi1}, Corollary~1.7,
and follows again as the left@-hand@-side is intrinsic of $Y$.
Proposition~1.3 is significant in view of the consequence:
\proclaim{Corollary 1.4}
$$\mu_{\Cal L}(Y)=(-1)^{\dim Y} \cwm(Y)_{\vee\Cal L}\quad.$$
\endproclaim
The class $\mu_{\Cal L}(Y)$ is the `$\mu$@-class' defined and studied in
\cite{Aluffi1}; it carries a notable amount of information about $X$,
with applications to duality and to the study of contacts of
hypersurfaces. Corollary~1.4 solves a puzzle left open in
\cite{Aluffi1} (p.~326): to define a class for arbitrary schemes,
specializing to $\mu_{\Cal L}(Y)$ for singular schemes of
hypersurfaces. It also clarifies the dependence of the $\mu$@-class on
the line bundle $\Cal L$: it follows from Corollary~1.4 that if $\Cal
L_1$, $\Cal L_2$ are line bundles, then
$$\mu_{\Cal L_2}(Y)=\mu_{\Cal L_1}(Y)_{\Cal L_1^\vee\otimes\Cal L_2}$$
(this does {\it not\/} follow formally from the expression for the
$\mu$@-class in terms of the Segre class of $Y$!) For applications of
Proposition~1.3 and Corollary~1.4, see Examples~3.4, 3.5.

The next fact we list also requires some notations. We now assume that
$X$ is a reduced hypersurface, over $\Bbb C$. The question is whether,
in this particularly `geometric' case, $\cwm(Y)$ can be recovered from
numerical invariants of $X$. The answer comes from \cite{P@-P}:
define a function $\mu: Y @>>> \Bbb Z$ by setting $\mu(y)=(-1)^{\dim
X}(\chi(y)-1)$, where $\chi(y)$ is the Euler characteristic of the
Milnor fiber of $X$ at $y$; $\mu$ is a constructible function on
$Y$, so we can apply to it MacPherson's transformation $\cmp$ (that is,
write $\mu$ as a linear combination of characteristic functions $1_Z$
for subvarieties $Z$ of $Y$, then replace each $1_Z$ in this
combination by $\cmp(Z)$).
\proclaim{Theorem 1.5}
$$\cwm(Y)=(-1)^{\dim Y}\cmp(\mu)\quad \text{in $A_*Y$.}$$
\endproclaim
Equivalently, write $\mu$ as a linear combination of local Euler
obstructions (also an ingredient in \cite{MacPherson}): $\mu=\sum\ell_i
\text{Eu}_{Y_i}$; then the content of Theorem~1.5 is that in this
situation the $Y_i$'s are precisely the supports of the components of
the normal cone of $Y$, and the numbers $\ell_i$ determined by $\mu$
agree (up to sign) with the multiplicities $m_i$ used to define
$\cwm(Y)$. Again, we would be very interested in extensions of this
result to more general $Y$: what numerical invariants of a space $X$
(not necessarily a hypersurface) determine the multiplicities of the
components of the normal cone of its singularity subscheme? Can these
multiplicities be computed for an arbitrary scheme $Y$, by a similar
`Milnor fiber' approach? Once more, note that the left@-hand@-side in
Theorem~1.5 is defined for arbitrary $Y$; to what extent can the
right@-hand@-side also be defined for arbitrary $Y$? We know of
several problems in enumerative geometry for which finding these
multiplicities is one of the main computational ingredients. For an
explicit computation (not directly related to enumerative geometry)
see Example~3.6.

Finally, it would be interesting to have results on the functoriality
of the class $\cwm(Y)$; little is known about the functoriality of the
ordinary Chern@-Mather class. Again, something can be said if $Y$ is the
singularity subscheme of a hypersurface $X$ (over an arbitrary
algebraically closed field of characteristic zero, and possibly
nonreduced). Let $Z$ be a nonsingular subvariety of $Y\subset X\subset
M$, and consider the blow@-up $\Til M$ of $M$ along $Z$:
$$\diagram
       & {Y'} \rto \dto & {X'} \rto \dto^{\rho} & {\Til M} \dto^{\pi} \\
Z \rto & Y   \rto      &  X   \rto      &     M
\enddiagram$$
Here $X'=\pi^{-1} X$ is the scheme@-theoretic inverse image of $X$, a
hypersurface of $\Til M$, and $Y'$ is the singularity subscheme of
$X'$.
\proclaim{Proposition 1.6} Assume $Z$ has codimension $d$ in $M$. Then
$$\rho_* \cwm(Y')=(-1)^{\dim X-\dim Y}\cwm(Y)-(d-1) \cwm(Z)\quad \text{in
$A_*X$.}$$
\endproclaim
Here of course $\cwm(Z)=c(TZ)\cap [Z]$, as $Z$ is nonsingular. Also note
that by assumption $X$ is singular along $Z$, hence $Y'$ contains the
exceptional divisor in $\Til M$. 
\vskip 6pt

Proofs of the statements made in this section are sketched in \S2,
with emphasis on Theorem~1.2, which relates the weighted Chern@-Mather
class of the singularity of a hypersurface with its {\it Milnor
class.\/}

 
\head \S2. The Milnor class of a hypersurface.\endhead

As is well known, for a compact complex hypersurface $X$ with isolated
singularities the sum of the Milnor numbers of the singularities
measures the difference between the topological Euler characteristic of
$X$ and that of a nonsingular hypersurface linearly equivalent to $X$
(if there is such a hypersurface). To my knowledge, the first who used
this fact to define and study a generalization of the Milnor number to
non@-isolated hypersurface singularities is Adam Parusi\`nski,
\cite{Parusi\`nski}. 

Now, the functoriality of Schwartz@-MacPherson's class implies that,
for a hypersurface $X$ as above, the Euler characteristic of $X$ equals
the degree of the (zero@-dimensional component of the) class $\cmp(X)$.
On the other hand, the Euler characteristic of a nonsingular
hypersurface linearly equivalent to $X$ equals the degree of the class
of the virtual tangent bundle of $X$ (that is, of $c_F(X)$ with
notations as in \S1). That is, Parusi\`nski's Milnor number equals (up
to a sign), the degree of the difference between the two classes:
$$\int\left(c_F(X)-\cmp(X)\right)\quad.$$
It is natural then to study the whole class $c_F(X)-\cmp(X)$; this (or
slight variations of it) has been named the {\it Milnor class\/} of $X$
by some authors (see \cite{B@-L@-S@-S}, \cite{P@-P}, \cite{Yokura}).

Note that nothing in the definition of the class $c_F(X)-\cmp(X)$
requires $X$ to be a hypersurface: both Schwartz@-MacPherson's and
Fulton's classes can be defined for arbitrary varieties. For reduced
compact complex local complete intersections, the Milnor class is
computed in homology in \cite{B@-L@-S@-S} in terms of vector fields on
$X$, an approach reminiscent of Schwartz's definition of $\cmp(X)$.

In fact the class makes sense for arbitrary schemes $X$ over any
algebraically closed field of characteristic~0, and naturally lives in
the Chow group $A_*Y$ of the singular locus of $X$. We would like to
pose the following question:

{\it ---To what extent is the Milnor class of $X$ determined by the
singular scheme $Y$ of $X$?} or, in more ambitious terms:

{\it ---Is there a natural definition of a class on an arbitrary scheme
$Y$, from which the Milnor class of $X$ can be computed if $Y$ is the
singular scheme of $X$?}

In view of the results collected in \S1, the situation is clear for
hypersurfaces. The singular locus of a hypersurface has a natural
scheme structure, given by the partial derivatives of local equations
of~$X$. Theorem~1.2 then asserts that (for arbitrary hypersurfaces $X$
over an algebraically closed field of characteristic~0, and writing
$\Cal L=\Cal O(X)_{|Y}$) 
$$c_F(X)-\cmp(X)=(-1)^{\dim X-\dim Y} c(\Cal
L)^{-1} \cap \cwm(Y)\quad \text{in $A_*X$:}$$
that is, {\it if two hypersurfaces have the same singular scheme $Y$
and their line bundles restrict to the same bundle on $Y$, then they
have the same Milnor class; and, further, this can be recovered from
the class $\cwm(Y)$, which can be defined for arbitrary schemes $Y$.}

Therefore, Theorem~1.2 answers the two questions posed above, for
hypersurfaces. To our knowledge, the questions are completely open for
more general schemes $X$. Milnor classes of local complete
intersections (for which the singular locus also carries a natural
scheme structure) have been studied in \cite{B@-B@-L@-S}, but from a
different viewpoint, which does not seem to address questions such as
the ones posed above.

Theorem~1.2 could be deduced easily from results in the existing
literature (particularly from \cite{P@-P} or \cite{Aluffi2}). However,
while the main result in \cite{Aluffi2} is at the right level of
generality, its proof is rather involved; on the other hand, the
approach in \cite{P@-P} is more streamlined, but it seems to rely on
the complex geometry of the situation, and to require the hypersurface
to be reduced. Also, the proof in \cite{P@-P} relies on the
sophisticated machinery of \cite{B@-M@-M}. The argument given below
works for possibly nonreduced hypersurfaces, over arbitrary
algebraically closed fields of characteristic~0, gives the formula in
rational equivalence, and only relies on the basic formalism of
Schwartz@-MacPherson's classes (as developed in \cite{Kennedy}). At
its core, however, is a multiplicity computation we learned from
\cite{P@-P}.

\demo{Proof of Theorem 1.2} We consider the blow@-up $\Til M
@>\pi>> M$ along $Y$, and let $\Cal X$, $\Cal Y$ be the pull@-back of
$X$ and the exceptional divisor, respectively. Note that $\Cal
Y\subset \Cal X$, so there is an effective Cartier divisor in $\Til M$
whose cycle equals $\Cal X-\Cal Y$; we will denote this divisor by
$\Cal X-\Cal Y$. Now let $p$ be a point of $X$. We have
$\pi^{-1}(p)\subset \Cal X-\Cal Y$, so it makes sense to consider the
Segre class of $\pi^{-1}(p)$ in $\Cal X-\Cal Y$.
\proclaim{Claim 2.1} Denoting degree by $\int$,
$$\int\frac{s(\pi^{-1}(p),\Cal X-\Cal Y)}{1+\Cal X-\Cal Y}=1$$
\endproclaim

A preliminary result is in order before we prove this claim. We have
$$\pi^{-1}(p) \hookrightarrow (\Cal X-\Cal Y) \hookrightarrow \Til M\quad,$$
where the second embedding is regular. We claim that
$$s(\pi^{-1}(p),\Cal X-\Cal Y)=c(N_{\Cal X-\Cal Y}\Til M)\cap
s(\pi^{-1}(p),\Til M)\quad.$$
Note that this is {\it not\/} automatic in this situation,
cf.~Example~4.2.8 in \cite{Fulton}. In our case, it will follow from
the following lemma:
\proclaim{Lemma 2.2} Let $D$, $E$ be hypersurfaces in a variety $V$.
Assume that $D-E$ is positive and has no components in common with
$E$. Then $s(E,D)=c(N_DV)\cap s(E,V)$.\endproclaim
\demo{Proof of the lemma} By the hypothesis and Lemma~4.2 in \cite{Fulton},
$$\align
s(E,D)&=s(E,E)+s(E\cap(D-E),D-E)=[E]+\frac{E\cdot(D-E)}{1+E}\\
&=\frac{([E]+E\cdot E)+E \cdot (D-E)}{1+E}=(1+D)\cap \frac{[E]}{1+E}\\
&=c(N_DV)\cap s(E,V)\quad.\qed
\endalign$$
\enddemo
\demo{Proof of Claim 2.1} We apply Lemma~2.2 to the normalized
blow@-up $V$ of $\Til M$ along $\pi^{-1}(p)$, with $E=$the exceptional
divisor, and $D=$the inverse image of $\Cal X-\Cal Y$. To see that the
hypotheses are satisfied, we have to show that every component of $E$
appears with the same multiplicity in $E$ and $D$. 

For this\footnote{This computation is essentially lifted from an
analogous computation in the proof of Proposition~2.2 in \cite{P@-P}.},
let $\gamma(t)$ be a germ of a nonsingular curve centered at the
general point of a component of $E$, let $\tilde\gamma(t)$ be the
composition to $M$, and let $F$ be a local equation
for $X$ at $p$; also, choose local parameters $x_1,\dots,x_n$
for $M$ at $p$. The ideal of $E$ is the pull@-back of
$(x_1,\dots,x_n)$ to $V$, so the multiplicity $m_E$ of the component
in $E$ equals the order of vanishing of the pull@-back
$x_i(t)=\tilde\gamma^*x_i$ of a generic local parameter. The multiplicity
$m_D$ in $D$ equals $m_{\Cal X}-m_{\Cal Y}$, where $m_{\Cal X}$,
$m_{\Cal Y}$ are respectively the multiplicities in the pull@-backs of
$\Cal X$, $\Cal Y$.

Now $m_{\Cal X}$ is the order of vanishing of 
$$\tilde\gamma^* F=F(x_1(t),\dots,x_n(t))\quad,$$
while $m_{\Cal Y}$ is the order of vanishing of the pull@-back of
$$\left(F,\frac{\partial F}{\partial x_1},\dots,\frac{\partial
F}{\partial x_n}\right)\quad,$$
that is, the order of vanishing of $\tilde\gamma^*\frac{\partial F}{\partial
x_i}$ for a generic local parameter $x_i$. Now taking the derivative
with respect to $t$ gives (by the chain rule!)
$$m_{\Cal X}-1=m_{\Cal Y}+m_E-1\quad,$$
from which the desired equality $m_E=m_D$ follows.

Applying Lemma~2.2, we get
$$s(E,D)=(1+\Cal X-\Cal Y)\cap s(E,V)\quad,$$
hence
$$s(\pi^{-1}(p),\Cal X-\Cal Y)=(1+\Cal X-\Cal Y)\cap s(\pi^{-1}(p),\Til M)$$
by the birational invariance of Segre classes (\cite{Fulton},
Proposition~4.2). From this,
$$\pi_*\frac{s(\pi^{-1},\Cal X-\Cal Y)}{1+\Cal X-\Cal Y}=s(p,M)=[p]\quad,$$
again by the birational invariance of Segre classes, and the claim
follows by taking degrees.\qed\enddemo

We are finally ready to prove Theorem~1.2. Identify $\Cal Y$ with
the projective normal cone of $Y$ in $M$, let $\Cal Y_i$ be the
reduced components of $\Cal Y$, and let $Y_i$ be their support in $Y$. Then
$$\left\{\aligned
\Cal X &=\Til X+\sum n_i\Cal Y_i\\
\Cal Y &=\sum m_i \Cal Y_i
\endaligned\right.$$
for suitable $m_i$, $n_i$. By Claim 2.1,
$$\align
1&=\int \frac{s(\pi^{-1}(p),\Cal X-\Cal Y)}{1+\Cal X-\Cal Y}\\
&= \int \frac{s(\pi^{-1}(p)\cap \Til X,\Til X)+\sum (n_i-m_i)
s(\pi^{-1}(p)\cap \Cal Y_i,\Cal Y_i)}{1+\Cal X-\Cal Y}\\
\intertext{by Lemma 4.2 in \cite{Fulton}}
&= \text{Eu}_X(p) +\sum (n_i-m_i) (-1)^{\dim M+1-\dim Y_i}
\text{Eu}_{Y_i}(p)
\endalign$$
using the formula for Euler obstructions due to
Gonzalez@-Sprinberg and Verdier, as computed in \cite{Kennedy}, Lemma~2
(as pointed out in \cite{Aluffi2}, \S1.3 and in \cite{P@-P}, \S3,
the divisor $\Cal X-\Cal Y$ can be embedded in $\P(T^*M)$, and $1+\Cal
X-\Cal Y$ is then the restriction of the class of the tautological
bundle in $\P(T^*M)$). Now, every relation between constructible
functions yields a relation for characteristic classes. Here, this
gives (using the formula for Mather's classes in \cite{Kennedy},
Lemma~1, going back to Claude Sabbah):
$$\align
\cmp(X) &=\cma(X)+\sum (n_i-m_i) (-1)^{\dim M+1-\dim Y_i}\cma(Y_i)\\
&=c(TM)\cap \pi_*\left(\frac{[\Til X]}{1+\Cal X-\Cal Y}+\sum
(n_i-m_i) \frac{[\Cal Y_i]}{1+\Cal X-\Cal Y}\right)\\
&=c(TM)\cap \pi_*\left(\frac{[\Cal X]}{1+\Cal X}-\frac 1{1+\Cal X}\sum m_i
\frac{[\Cal Y_i]}{1+\Cal X-\Cal Y}\right)\\
&=c_F(X)+c(\Cal L)^{-1}\cap \sum m_i (-1)^{\dim M-\dim Y_i} \cma(Y_i)\\
&=c_F(X)+(-1)^{\dim M-\dim Y} c(\Cal L)^{-1}\cap \cwm(Y)\endalign$$
which is the desired formula.\qed\enddemo

As observed in the proof, $\Cal X-\Cal Y$ can be naturally embedded in
$\P(T^*M)$. The content of Claim 2.1 is that $\Cal X-\Cal Y$ gives
then the {\it characteristic cycle\/} of $X$ (corresponding to the
characteristic function $1_X$ of $X$ in $M$).

The other statements in \S1 now follow easily, either by comparing the
expression for $\cmp$ with the expressions in \cite{Aluffi2} and
\cite{P@-P}, or by direct manipulations that can be extracted from those
sources. The argument given here re@-proves Theorem~I.3 in
\cite{Aluffi2}/Theorem~3.1 in \cite{P@-P}; and for example, \S1 in
\cite{Aluffi2} shows how to go directly from this form of the result
to expressions in terms of Segre classes or $\mu$@-classes (thus
proving Proposition~1.3, Corollary~1.4). 

The details are left to the reader. Theorem~1.5 is our reading of
Theorem~2.3 (iii) from \cite{P@-P}. The blow@-up formula of
Proposition~1.6 follows from Proposition~IV.2 in \cite{Aluffi2}.


\head \S3. Examples and applications.\endhead

Normal cones behave well with respect to proper finite maps and with
respect to flat maps, cf.~Proposition~4.2 in \cite{Fulton}. For
example, assume that $Y$ is irreducible, and $\Til M @>\pi>> M$ is a
surjective birational map on the ambient space. Then there is an
induced surjective map from the cone of $\pi^{-1} Y$ to the cone of
$Y$. This can be used to obtain the data $\{(Y_i,m_i)\}$ of \S1, for
example by suitably blowing up an ambient space; this can lead to
direct computations of weighted Chern@-Mather classes.

\example{Example 3.1} Suppose $Y$ consists of a curve $C$, with an
embedded multiple planar point at a point $p$. More precisely, assume
$C$, $Y$ have local ideals respectively $\Cal I_C$, $\Cal I_C\cdot
(x,y)^m$, $m\ge1$, near $p$ in a nonsingular ambient surface $S$ with
local parameters $x$, $y$. Also, assume that $C$ has multiplicity $r$
at $p$. Then
$$\cwm(Y)=\cma(C)-(m+r)[p]\quad.$$
Indeed, blow@-up $S$ at $p$; the total transform of $Y$ consists of the
proper transform of $C$, and of $(m+r)$ times the exceptional divisor.
Therefore, the normal cone of $Y$ contains a component with multiplicity
$m+r$ over $p$. (But note there is no such component if $m=0$.)

For example, take $Y$ to be the union of two lines $\ell_1$, $\ell_2$
in $\P^2$, with an embedded planar point at the intersection
$p=\ell_1\cap\ell_2$; then $\cwm(Y)=[\ell_1]+[\ell_2]+[p]$. If the
embedded point is on one of the lines, but not at $p$, then
$\cwm(Y)=[\ell_1]+[\ell_2]+2[p]$. If each line comes with multiplicity
$r$, and the embedded point is at $p$, then the class is
$$r \cma(\ell_1)+r \cma(\ell_2)-(1+2 r)[p]=r[\ell_1]+r[\ell_2]+
(2r-1)[p]\quad.$$
\endexample

\example{Example 3.2} Example~3.1 can be easily generalized to the
situation in which $Y$ is a subscheme of a given ambient space $M$,
and the residual to a Cartier divisor $D$ in $Y$ is a known scheme
$Y'$. Then $\cwm(Y)$ can be written in terms of $\cma(D)$, $\cwm(Y')$,
and the multiplicity of $D$ along the distinguished components of
$Y'$; details are left to the reader. A very different expression can
be obtained if $Y'$ is the singularity subscheme of a hypersurface $X$
in $M$, and $D$ is the $r$@-th multiple of~$X$ ($r\ge 0$).

\proclaim{Claim 3.1} Let $\Cal L=\Cal O(X)_{|Y'}$. Then
$$\cwm(Y)=r\,c_F(X)+(-1)^{\dim X-\dim Y'}\frac{c(\Cal L^{\otimes(r+1)})}
{c(\Cal L)}\cap \cwm(Y')\quad.$$
\endproclaim
The proof is an easy application of the results in \S1, and is also
left to the reader. 

To contrast the two approaches, take again the example of the union of
two lines $\ell_1$, $\ell_2$ in $\P^2$, each coming with
multiplicity~$r$, with an embedded planar point at the intersection.
Since the planar point is the singularity subscheme of the union of
two (simple) lines, Claim~3.1 computes the weighted Chern@-Mather
class of this scheme as
$$\align
r\, c_F(X)+&(-1)^{\dim X-\dim Y'}\frac{c(\Cal L^{\otimes(r+1)})}
{c(\Cal L)}\cap \cwm(Y')
=r \frac{c(T\P^2)}{c(\Cal O_{\P^2}(2))}\cap ([\ell_1]+[\ell_2]) -[p]\\
&=r([\ell_1]+[p])+r([\ell_2]+[p])-[p]
\endalign$$
with the same result as before, but by a very different route.

It would be useful to have formulas such as Claim~3.1, but with less
stringent hypotheses on $X$.\endexample

\example{Example 3.3}
If $X=X_1\cup \cdots \cup X_r$ is a divisor with normal crossings,
with all $X_i$  supported on nonsingular hypersurfaces
$(X_i)_{\text{red}}$, and $Y$ is its singularity subscheme, then
$$\cwm(Y)=\pm c(TM)\cap\left(1-\frac{1+[X]}{(1+(X_1)_{\text{red}})
\cdots (1+(X_r)_{\text{red}})}\right)\cap [M]\quad,$$
taking the sign $+$, resp.~$-$ according to whether $X$ is reduced or
not. The expression is interpreted by expanding it, which leaves a
class naturally supported on $Y$; it follows from Proposition~1.3
and \cite{Aluffi2}, \S2.2 (Lemma~II.2 in \cite{Aluffi2} computes the
Segre class if $X$ is reduced, and the computation in the proof of
Lemma~II.1 is used to cover the non@-reduced case).\endexample

\example{Example 3.4}
What do we learn about hypersurfaces by studying their Milnor classes?

As shown in \cite{Aluffi1}, the $\mu$@-class of a hypersurface $X$ packs
a good amount of information about $X$: for example, the multiplicity of
$X$ as a point of the discriminant of a linear system and the dimension
of this discriminant can be recovered very easily from the $\mu$@-class
(hence from the Milnor class). In the classical language, the
$\mu$@-classes of hyperplane sections of an embedded nonsingular
projective variety $M$ gives a localized analog of the {\it ranks\/}
of $M$, and provides a natural tool to study projective duality.

In a different direction, the good behavior of the $\mu$@-class can be
used to put restrictions on the possible singularities of a
hypersurface in a given ambient space. Several examples of this
phenomenon are illustrated in \cite{Aluffi1}, \S3,
where the main tool was the observation that if the singularity subscheme
$Y$ of a hypersurface $X$ is {\it nonsingular,\/} then
$$\mu_{\Cal L}(Y)=c(T^*Y\otimes\Cal L)\cap [Y]\quad.$$
Now, Corollary~1.4 from \S1:
$$\mu_{\Cal L}(Y)=(-1)^{\dim Y} \cwm(Y)_{\vee \Cal L}$$
is a substantial upgrade of this formula, and this allows us to extend
some of those results.

\proclaim{Claim 3.2} If two smooth hypersurfaces of degree $d_1$, $d_2$ in
projective space are tangent along a positive dimensional set, then
$d_1=d_2$.\endproclaim

More generally, if two smooth hypersurfaces $M_1$, $M_2$ of a variety
$V$ are tangent along an irreducible (for simplicity) set $Z$, and
$\dim Z>0$, then we claim that 
$$r\,M_1\cdot [Z]=r\,M_2\cdot [Z]\quad,$$
where $r$ is the order of tangency of $M_1$ and $M_2$ (for example,
$r=1$ if $M_1$, $M_2$ have simple contact). This is essentially
Proposition~IV.7 in \cite{Aluffi2}, with all hypotheses on the contact
locus (except the positive dimensionality) removed. The stronger
statement given above follows from the results in \S1. Indeed, in the
situation of the statement, let $X=M_1\cap M_2$; then $X$ is a
hypersurface in two distinct ways: with respect to $\Cal L_2=\Cal
O(M_2)_{|M_1}$ in $M_1$, and with respect to $\Cal L_1=\Cal
O(M_1)_{|M_2}$ in $M_2$. The contact locus is $Y=\text{Sing}\,X$ (with
the scheme structure specified in \S1), and $[Y]=r[Z]$. By Theorem~1.2
$$ c(\Cal L_2)^{-1}\cap \cwm(Y)=c(\Cal L_1)^{-1}\cap \cwm(Y)\quad,$$
implying
$$c_1(\Cal L_1)\cap [Y]=c_1(\Cal L_2)\cap [Y]\quad,$$
which is the statement.\endexample

\example{Example 3.5} We say that a hypersurface $X$ of a nonsingular
variety $M$ is (analytically) `homogeneous at $p$' if the equation of
$X$ is homogeneous for some choice of system of parameters in the
completion of the local ring for $M$ at $p$. We are going to consider
degree@-$d$ hypersurfaces $X$ in $\P^n$, whose singular scheme $Y$ has
a connected component supported on a nonsingular curve $C$ of genus
$g$ and degree $r$; we assume that $Y$ has the reduced structure at
all but finitely many points $q_1,\dots,q_s$, and that $X$ is
homogeneous at each of the $q_i$. In particular, $X$ has multiplicity
2 at all other points of $C$; we let $m_i$ be the multiplicity of $X$
at $q_i$.

How constrained is this situation? Examples~3.4---3.6 in \cite{Aluffi1}
deal with the case in which the singular scheme is reduced, that is,
there are no points `$q_i$' as above. This situation is then very
rigid: for example, one sees that only quadrics can have singular
scheme equal to a line, and no hypersurface in projective space
can have singular scheme equal to a twisted cubic (cf.~p.~347 in
\cite{Aluffi1}).

The natural expectation would be that letting the singular scheme 
be nonreduced should allow many more examples. For instance, cones
over nodal plane curves give examples of hypersurfaces in $\P^3$ of
arbitrary degree $\ge 2$ and singular scheme generically reduced, but
with an embedded homogeneous point (at the vertex). However, the
results in this paper show that the situation is still quite rigid:

\proclaim{Claim 3.3} Under the hypotheses detailed above, $(n-1)$ must
divide $4(g+r-1)$. In fact, necessarily
$$(n-1)\left((d-2)r-\sum(m_i-2)\right)=4(g+r-1)\quad.$$
\endproclaim

For example, twisted cubics can support singularity subschemes as
above only in dimensions $n=3, 5, 9$, regardless of the number of
embedded points allowed on them. (We do not know if such examples do
exist.) The only situation in unconstrained dimension is for
$g+r-1=0$, that is, $g=0$ and $r=1$: lines are the only nonsingular
curves in projective space which may support a generically reduced
singularity subscheme in all dimensions (under the local homogeneity
assumption). Further, if $Y$ is supported on a line and only has {\it
one\/} embedded homogeneous point, then the formula implies that the
multiplicity of $X$ at this point is $d$; therefore, $X$ is
necessarily a cone in this case.

For $\sum(m_i-2)=0$, the formula in Claim~3.3 recovers the formula at
p.~347 of \cite{Aluffi1} (that is, the reduced case). For $n=2$, the
hypotheses imply that $X$ is a plane curve consisting of a double
component $C$ and a residual curve of degree $(d-2r)$; the formula
then follows from the genus formula and B\'ezout's theorem. In higher
dimensions, the following argument is the only proof we know.

\demo{Proof of the claim} We compute directly the weighted
Chern@-Mather class of $Y$ and the Segre class $s(Y,\P^n)$.
Proposition~1.3 gives a relation between these two classes, and the
formula follows by taking degrees.

Explicitly, blow@-up $\P^n$ at the `special' points $q_1,\dots, q_s$,
and then along the proper transform of the curve $C$. The homogeneity
hypothesis implies that the (scheme@-theoretic) inverse image of $Y$
in the top blow@-up is a Cartier divisor, with a component of
multiplicity 1 dominating $C$, and $s$ components with multiplicity
$(m_1-1),\dots,(m_s-1)$ dominating the $q_i$'s. The Segre class of $Y$
in $\P^n$ is then computed by using the birational invariance of Segre
classes, and we get{\eightpoint
$$i_*s(Y,\P^n)=r[\P^1]+\left(s(n-1)+2-2g-r(n+1)+\sum_i
((m_i-1)^n-n(m_i-1))\right) [\P^0]$$}
\noindent (where $i:Y\hookrightarrow\P^n$ is the inclusion).

On the other hand, the component dominating $q_i$ maps to a
corresponding component of the projective normal cone to $Y$ in
$\P^n$; computing differentials, we see that this map has degree
$(m_i-1)^{n-1}-1$. This allows us to compute the weighted
Chern@-Mather class of $Y$:
$$\cwm(Y) =\cma(C)-\sum_i \left((m_i-1)^{n-1}-1\right)(m_i-1)
\cma(q_i) \quad,$$
from which
$$i_*\cwm(Y)=r[\P^1]+\left(2-2g-\sum_i ((m_i-1)^n-(m_i-1))\right)[\P^0]$$

Now let $h$ denote the hyperplane class in $\P^n$. The expression for
the Segre class gives
{\eightpoint $$\multline
i_*c(T^*\P^n\otimes\Cal O(d))\cap s(Y,\P^n) = i_*\frac{(1+(d-1)
h)^{n+1}}{1+d h}\cap s(Y,\P^n)\\
=r[\P^1]+\left((s+r d-2 r)(n-1)+2-2g-4r+r d+\sum_i
((m_i-1)^n-n(m_i-1))\right) [\P^0]
\endmultline$$}
and therefore
$$\multline
i_*(-1)^{\dim Y}\left(c(T^*M\otimes \Cal L)\cap s(Y,M)\right)_{\vee\Cal L}
=r[\P^1]\\
+\left((2r-dr-s)(n-1)-2+2g+4r-\sum_i ((m_i-1)^n-n(m_i-1))\right) [\P^0].
\endmultline$$

By Proposition 1.3, this class must equal $i_*\cwm(Y)$. Equating the
two expressions gives the formula in the
statement.\qed\enddemo\endexample

\example{Example 3.6} Finally, we give an example of the use of
weighted Chern@-Mather classes in the computation of the
multiplicities of components of a normal cone. Such multiplicities are
important for enumerative applications, and it would be very useful to
develop tools to compute them. For singularity subschemes of
hypersurfaces, the connection between weighted Chern@-Mather classes
and Milnor classes often lets us recover these multiplicities from
computations of MacPherson's classes and local Euler obstructions. It
would be interesting to extend such techniques to more general
schemes.

Let $D$ be the hypersurface of $\P^9$ parametrizing {\it singular\/}
plane cubics, and let $Y$ be its singularity subscheme. The following
picture represents the natural stratification of $D$ (with arrows
denoting specialization):
$$\epsffile{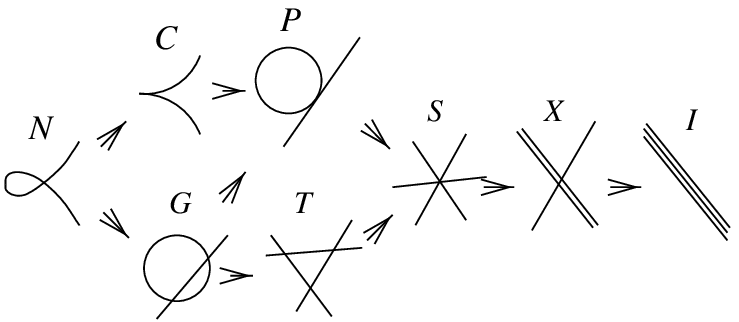}$$
The scheme $Y$ is supported on the union of the closures $\overline
C$, $\overline G$ of the loci parametrizing cuspidal cubics and
binodal cubics. What are the multiplicities of the components of the
normal cone of $Y$ in $\P^9$? The point here is that we can compute
$\cwm(Y)$ without knowing these multiplicities:

\proclaim{Claim 3.4} Denote by $i$ the inclusion of $Y$ in $\P^9$. Then
$$i_*\cwm(Y)= 69 [\P^7]+120 [\P^6]+210 [\P^5]+252 [\P^4]+210 [\P^3] +
120[\P^2]+45 [\P^1]+10 [\P^0].$$
\endproclaim
\demo{Proof} This follows from Theorem~1.2 and the computations of
characteristic classes for $D$ in \S4 of \cite{Aluffi4}.\qed\enddemo
Now the task is to find the coefficients expressing the weighted
Chern@-Mather class of $Y$ as a combination of the Chern@-Mather
classes of the loci $C$, $G$, etc. We first find the constructible
function $\nu$ corresponding to $\cwm(Y)$ under MacPherson's 
transformation. For this, we use the result of the computation from
\cite{Aluffi4} of Chern@-Schwartz@-MacPherson's classes of the strata
of $D$. Writing
$\cwm(Y)=\cmp(\nu)=\nu(C)\cdot \cmp(1_C)+\nu(G)\cdot \cmp(1_G)+\dots$
and solving the resulting system of linear equations, we find
$$\nu(C)=2;\,\, \nu(G)=1;\,\, \nu(P)=0;\,\, \nu(T)=1;\,\,
\nu(S)=3; \,\, \nu(X)=1;\,\, \nu(I)=1.$$
(The paragraph preceding the statement of Theorem~1.5 gives a
geometric interpretation of $\mu=-\nu$.) As pointed out in the
discussion following Theorem~1.5, to find the multiplicities we now
need to express this constructible function as a combination of local
Euler obstructions of the strata. These are easy to compute in
codimension~one, and we proceed to the computation of the multiplicities
for the components dominating the loci $\overline C$, $\overline G$,
$\overline P$, $\overline T$. For these loci, we only need to observe
that $\overline C$, $\overline G$ are nonsingular along $P$, and
$\overline G$ has multiplicity~3 along $T$ (these follows from easy
local computations). As the local Euler obstruction agrees with the
multiplicity in codimension~one, this gives{\eightpoint
$$\text{Eu}_C=\left\{\aligned
&\cdots\\
0\quad &T\\
1\quad &P\\
0\quad &G\\
1\quad &C
\endaligned\right.\quad,\quad
\text{Eu}_G=\left\{\aligned
&\cdots\\
3\quad &T\\
1\quad &P\\
1\quad &G\\
0\quad &C
\endaligned\right.\quad,$$}
\noindent where we indicate the value of the function at the general
point of the listed locus. Therefore
$$\nu=2\text{Eu}_C+\text{Eu}_G-3\text{Eu}_P-2\text{Eu}_T+\dots$$
from which we read that the multiplicities of the components of the
normal cone are: 2 over $C$, $1$ over $G$, $3$ over $P$, $2$ over $T$.

Finding the multiplicities over the remaining three loci $S$, $X$, $I$
requires computing the local Euler obstructions for all the strata of
$D$. We leave this to the motivated reader.\endexample



\enddocument